\documentclass[11pt]{article}
\usepackage[margin=2.5cm]{geometry}
\usepackage{amsmath,amsthm,amsfonts,amssymb}
\usepackage[mathscr]{euscript}
\usepackage{mathtools}
\usepackage[utf8]{inputenc}
\usepackage[T1]{fontenc}

\usepackage{graphicx,latexsym}
\usepackage{lscape}
\usepackage{inputenc}
\usepackage{multicol}
\usepackage{textcomp}
\usepackage{float}
\usepackage{graphics}

\newtheorem{thm}{Theorem}[section]
 \newtheorem{cor}[thm]{Corollary}
 \newtheorem{lem}[thm]{Lemma}
 \newtheorem{prop}[thm]{Proposition}
 \newtheorem{obsr}[thm]{Observation}
\newtheorem{ex}[thm]{Example}
 \newtheorem{defn}[thm]{Definition}

\newtheorem{rem}[thm]{Remark}

\numberwithin{equation}{section}

\newtheorem{discu}{Discussion:}

\newtheorem{conje}{Conjecture:}

\begin{document}
\date{}
\title{  $k$ -Strong Shortest Path Union Cover  for Certain  Graphs and Networks}
\author{\begin{tabular}{rcl}
           Antony Xavier$^1$, Santiagu Theresal$^1$, Deepa Mathew$^1$, S. Arul Amirtha Raja $^2$\\
        \end{tabular}\\
        \begin{tabular}{c}
        $^1$ Department of Mathematics, Loyola College, Chennai 600 034, India \\    Affiliated to University of Madras\\
         $^2$ Department of Mathematics, SRM Institute of Science and Technology, Chennai-89\\
         Email: santhia.teresa@gmail.com. \\
        \end{tabular}}
\maketitle
\vspace{-0.5 cm}
\begin{abstract}
The   $k$-distance  strong  shortest path union cover of a graph is the minimum cardinality among all strong shortest path union cover at distance $k$ of  $G$ where $1\leq k \leq d$. In this paper  we determine  the $2$-strong shortest path union cover  for certain  graphs, also we prove that the $k$-strong shortest path union cover problem  in general is NP-complete. 
\end{abstract}
\section{Introduction}\label{sec:Introduction}
\hspace*{0.5cm}    
For a graph $G$ = $(V,E)$, where $V$  and $E$ denote the vertex and the edge sets of a graph $G$, respectively. For vertices $u$ and $v$ in a connected graph G, the distance $d(u, v)$ is the length of a shortest $u - v$ path in $G$. The degree of a vertex $v$  in $G$ is the number of edges connecting it. The maximum distance between a vertex to all other vertices in a graph $G$ is considered as the eccentricity of the vertex. The diameter of a graph is the maximum eccentricity of any vertex in the graph. That is,  the greatest distance between any pair of vertices.  For more information on premlinaries  we refer the reader to \cite{IGT81}. Shortest path union cover problem was introduced by Peter Boothe et al \cite{PZA}. If $G$ is a graph, then a set $S$ of its vertices is called shortest path union cover, if the shortest  paths that start at the vertices of $S$ cover all the edges of $G$. The shortest  path union cover problem is to find  a shortest path union cover of minimum cardinality \cite{PZA}. Strong shortest path union cover problem is a new concept introduced by Xavier et all \cite{ASD}. Xavier et all \cite{ASD, ASDS, SXM}  proved that strong shortest path union cover problem is $NP$-complete  and obtained the results for certain graphs,  networks, Sierpi\'{n}skiGraphs and product graphs.\\
This manuscript is organized as follows.  In section $2$, the basic definitions and general results on $k$-strong shortest path union cover problem for certain graphs are obtained. In Section $3$, we discuss the Complexity of the $k$-strong shortest path union cover.  In section $4$, we determine the $2$-strong shortest path union cover problem for  certain graphs.  In section $5$, we discuss about the $2$-strong shortest path union cover problem for various networks and Sierpi\'{n}ski Graphs. 
\section{General  Results }
\begin{figure}[h!] 
	\centering
	\includegraphics[scale=0.6]{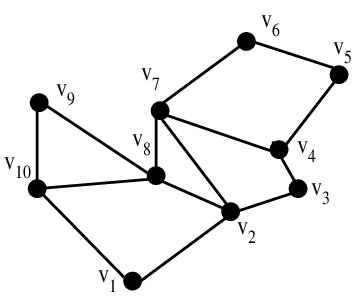}
	\caption {Graph $G$}
	\label{fig:Defin}
\end{figure}

\begin{defn}
Let $G$ be a graph, then a set $S$ of its vertices is called $k$-shortest path union cover, if the shortest  paths of length at most $k$ that start at the vertices of $S$ cover all the edges of $G$. The minimum cardinality of $k$-shortest path union cover  is denoted as  $SPC_{kU} (G)$.
\end{defn}	
\begin{defn}
	For a fixed vertex $u$, $\widetilde{P}_{k}(u,v)$ be the set of all edges of  a fixed shortest path between $u$ and $v$  where $v\in V(G)$ and $d(u,v)\leq k$.  Let  $\widetilde{P}_{k}(u)=\bigcup _{v \in V(G)}   \widetilde{P}_{k}(u,v)$. For any vertex set $S \subset V(G)$, define $\widetilde{P}_{k}(S)= \bigcup_{u \in S}   \widetilde{P}_{k}(u)$. If $\widetilde{P}_{k}(S) = E(G)$, then $S$ is called the $k$-strong shortest path union cover of $G$.The minimum cardinality of $k$-strong shortest path union cover is denoted as  $SSPC_{kU} (G)$.
\end{defn}	
\begin{ex}
In Figure \ref{fig:Defin}, $S=\{v_1, v_5,  v_{7}\}$ forms $2$-strong shortest path union cover and $T=\{v_1,  v_{7}\}$ forms $2$-shortest path union cover.
\end{ex}

\begin{thm}
	For any graph $G$, $\gamma_{k}(G)\leq SSPC_{kU}(G)$.
	\begin{proof}
		Every $k$-strong shortest path union cover of $G$ is a $k$-dominating set of $G$. Hence $\gamma_{k}(G)\leq SSPC_{kU}(G)$.			
	\end{proof}		
\end{thm}
\begin{obsr}
		For any  graph $G$ with  $\delta$$\geq$ $1$, $SSPC_{2U}(G)\leq\dfrac{n}{2}$
\end{obsr}
\begin{thm}
	For any  connected graph $G$, $SPC_{U}(G)\leq SPC_{kU}(G)\leq SSPC_{kU}(G)$.
	\begin{proof}
		Every $k$-shortest path union cover is a shortest path union cover of $G$.\\ Therefore,  $SPC_{U}(G)\leq SPC_{kU}(G)$. Every $k$-strong shortest path union cover is a $k$-shortest path union cover. Therefore, $SPC_{kU}(G)\leq SSPC_{kU}(G)$. \\
		Hence $SPC_{U}(G)\leq SPC_{kU}(G)\leq SSPC_{kU}(G)$.
	\end{proof}			
\end{thm}

\begin{lem}
	For any connected graph $G$ of order $n$, $SSPC_{U}(G) \leq SSPC_{kU}(G) \leq n-1$.
	\begin{proof}
		Every $k$-strong shortest path union cover of $G$ is a strong shortest path union cover of $G$. Also $V(G) /\{u\}$ where $u\in V(G)$ is a $k$-strong shortest path union cover of $G$. Therefore $SSPC_{kU}(G) \leq n-1$. Thus 	 $SSPC_{U}(G) \leq SSPC_{kU}(G) \leq n-1$.
	\end{proof}		
\end{lem}

\begin{rem}\label{re1}
	If $G$ is any connected graph with $\delta$$\geq$ $2$, then $SSPC_{kU}(G)$$\geq$$2$.
\end{rem}

\begin{rem}\label{re3}
If $G$ is a graph with no pendent vertices then $SSPC_{kU}(G) \geq 2$.	
\end{rem}

\begin{figure}[H] 
	\centering
	\includegraphics[scale=0.6]{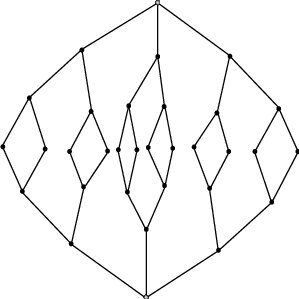}
	\caption {Graph $G$}
	\label{fig:Diam}
\end{figure}

\begin{thm}\label{bou}
	For a graph $G$, with   $SPC_{kU}(G)\geq \frac{|E|(\Delta-2)}{\Delta ((\Delta-1)^{k}-1)}$. 
	\begin{proof}	
		Let $S\subseteq V(G)$ be the $k$-shortest path union cover of $G$. This implies that every edge of $G$ lies on a path with distance at most $k$ from the vertices in $S$. Consider a vertex $u\in S$. All  possible paths of length at most $k$ with one end $u$ can cover at most $\Delta+\Delta(\Delta-1)+\Delta(\Delta-1)^2+\ldots+\Delta(\Delta-1)^{k-1}$ edges. That is the vertex $u$ can cover at most $\Delta \frac{(\Delta-1)^{k}-1}{(\Delta-2)}$ edges. Therefore $|E|\leq|S|\Delta \frac{(\Delta-1)^{k}-1}{(\Delta-2)}$. Hence $|S|\geq\frac{|E|(\Delta-2)}{\Delta ((\Delta-1)^{k}-1)} $.
	\end{proof}	
\end{thm}
\begin{rem}
The bound in Theorem \ref{bou}	is sharp for the graph given in Figure \ref{fig:Diam}, when $k=3$.
\end{rem}

\begin{figure}[H] 
	\centering
	\includegraphics[scale=0.6]{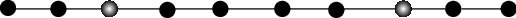}
	\caption {Graph $G$}
	\label{fig:Diam1}
\end{figure}

\begin{figure}[H] 
	\centering
	\includegraphics[scale=0.6]{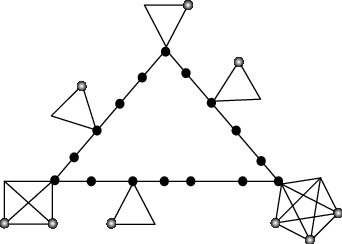}
	\caption {Graph $G$}
	\label{fig:Clique}
\end{figure}

\begin{thm}
	If $G$ is a nontrivial connected graph of order $n$ and diameter $d$, then $SSPC_{kU} (G)  \leq n - k + 1$ where $k\leq d$.
	\begin{proof}
	Since $diam(G)=d$, there exists a $u-v$ geodesic of length $k$. Let $u= v_0, v_1, v_2, \dots v_k= v$ be the $u-v$ geodesic. Let $S$ = $V(G)$ -\{$v_1, v_2, \dots v_{k-1}$\}. Then all the edges of $G$ are  coverd by the vertices in $S$. Hence $SSPC_{kU} (G)$ $\leq$ $|S|$= $n - k + 1$.
	\end{proof}	
\end{thm}
\begin{thm}\label{bou1}
For a graph $G(V,E)$ with diameter $d \geq 2$, the $k$-strong 
shortest path union cover $SSPC_{kU} (G) \leq n-(d+1)+(\frac{d+1}{2k+1})$.
	\begin{proof}
	Consider a graph $G(V(G),E(G))$ with diameter $d \geq 2$. Let $u$ and $v$ be the vertices of $G$ for which $d(u,v) =d$. Assume that $S\subseteq V(G)$ forms a $k$-strong shortest path union cover set for $G$. 
	Let $u = v_{k}, v_{3k}, \dots v_{(d+1)}= v$  be a $u-v$ 
	path of length at  distance $k$. Let $S = V(G) -\{v_{k}, v_{3k}, \dots v_{(d+1)}\}$. Thus  $SSPC_{kU} (G) \leq |S| \leq |V(G)|-(d+1)+(\frac{d+1}{2k+1})$.
		\end{proof}	
\end{thm}
\begin{rem}
	The bound in Theorem \ref{bou1}	is sharp for the graph given in Figure \ref{fig:Diam1}, when $k=2$.
\end{rem}
\begin{thm} \label{cli}
	$G$ has cliques $K_{n_{1}}$, $\dots$, $K_{n_{n}}$ each cliques has $n_{i}'$ simplicial vertices then $SSPC_{kU}(G)$ $\geq$$\sum$ $(n_{i}' - 1)$, where $n_{i}'$$\geq$$2$.
	\begin{proof}
		Consider the clique $K_{n_{i}}$	with simplicial vertices $n_{i}'$. Let $K_{n_{i}}'$ be the graph spanned by $n_{i}'$ simplicial vertices. It is straightforward to see that any geodesic of length at most $k$ through the non simplicial vertices of $K_{n_{i}}$ will not cover any edge of $K_{n_{i}}'$. To cover the edges of $K_{n_{i}}'$,  requires $n_{i}^{'}-1$ simplicial vertices of $K_{n_{i}}'$. Therefore $SSPC_{kU}(G)$ $\geq$$\sum$ $(n_{i}' - 1)$.
	\end{proof}		
\end{thm}
\begin{rem}
 The bound in Theorem \ref{cli}	 is sharp for the graph given in Figure 	\ref{fig:Clique} .
\end{rem}
\section{Complexity Results}
\begin{figure}[h!] 
	\centering
	\includegraphics[scale=0.32]{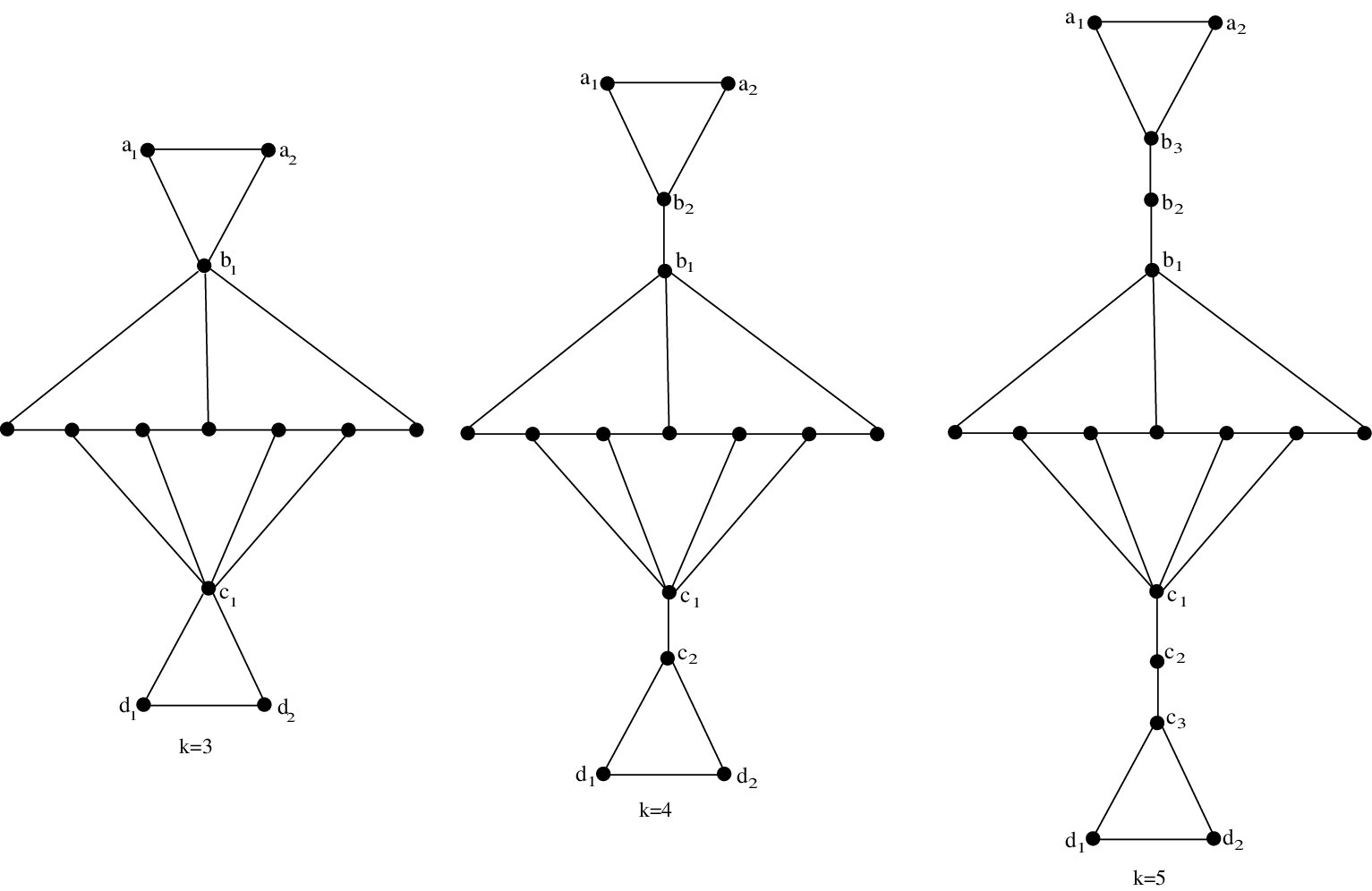}
	\caption {NP-complete illustration for a general graph}
	\label{fig:knpc2}
\end{figure}

The proof for the NP-completeness of the $k$-strong shortest path union cover problem for general graphs can be reduced from the
 vertex cover problem which is already proved to NPcomplete \cite{LH83}.
A vertex cover in an undirected graph $G$ = $(V, E)$ is a subset of
vertices $V$$^{'}$$\subseteq$ $V$ such that every edge in $G$ has at least one endpoint in $V^{'}$. 
	\begin{figure}[h!] 
	\centering
	\includegraphics[scale=0.3]{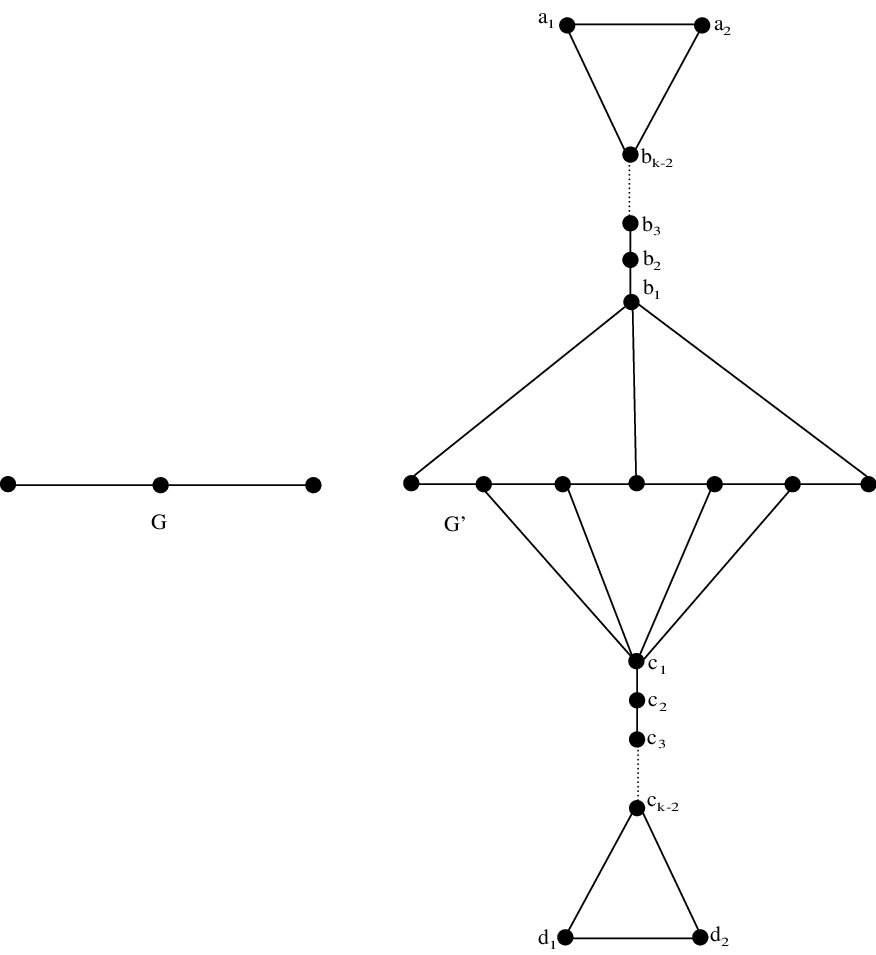}
	\caption {NP-complete illustration for a general graph}
	\label{fig:knpc}
\end{figure}
\begin{thm}
	The $k$-strong shortest path union cover problem is $NP$ – Complete.
	\begin{proof}
	There are two cases to prove that  the  $k$-strong shortest path union cover problem is $NP$-Complete, when  $k\geq3$ and $k=2$.\\
		Case (i): When $k\geq3$\\
		Given a graph $G$ = $(V, E)$.  We construct a graph $G' = (V', E')$  from $G$ by subdividing every edge in $E(G)$ by a path of length $3$.  All the new vertices that are created in this step are called   path vertices. Now create  two new $K_{3}$  and introduce two new vertices $b_{1}$ and $c_{1}$ such that all the vertices of  $V(G)$ in $G^{'}$ will be joined with $b_{1}$ and all the path vertices in $G^{'}$ will be joined with vertex $c_{1}$. Join any one vertex from anyone $K_{3}$ to $b_{1}$  with the path $b_{1}b_{2} \dots b_{k-2}$ and anyone vertex of other $K_{3}$  to $c_{1}$ with the path $c_{1}c_{2} \dots c_{k-2}$ of length $k-3$ as shown in the Figure 	\ref{fig:knpc}. Thus the  graph $G'$ is obtained.\\
		\hspace*{1.0cm}  
		Let $S$ be the vertex cover of $G$. Let $A=\{x,y /x$ = $a_{1}$ or $a_{2}$\\
\hspace*{8.8cm} 	    y = $d_{1}$ or $d_{2}$\}.\\
It is straightforward to see that $S\cup A$ is a $k$-strong shortest path union cover of $G'$ and $|S\cup A|=|S|+2$.\\
	From the construction of $G'$, any $k$-strong shortest path union cover of $G$ contains one vertex each from  \{$a_{1}, a_{2}$\} and \{$d_{1}, d_{2}$\} to cover the edges $a_{1}a_{2}$ and $d_{1}d_{2}$.\\	
	Conversely let $T$ be a $k$-strong shortest path union cover of  $G'$. Let $T' = T\setminus  \{b_{1},b_{2},\dots, b_{k-2}\} \cup \{c_{1},c_{2},\dots, c_{k-2}\}$. Note that $|T'|\leq |T|$. Construct $T''$ as follows.
	If  $T''$, contains 	$a_{1}$ and  $a_{2}$, then consider either $a_{1}$ or $a_{2}$ in $T''$. If $T'$ contains only one vertex from  \{$a_{1}, a_{2}$\}, then consider the same vertex in  $T''$. Similarly, If  $T'$, contains 	$d_{1}$ and  $d_{2}$, then consider either $d_{1}$ or $d_{2}$ in $T''$. If $T'$ contains only one vertex from  \{$d_{1}, d_{2}$\}, then consider the same vertex in  $T''$. It is easy to verify that  $|T''|\leq |T'|$.\\ 
	Now if $T''$ contains a path vertex $x'$, then replace $x'$ with the adjacent vertex $x \in V(G)$. Repeat this replacement procedure  until there are no more path vertices and the resulting set  be  $T^{'''}$. Also $|T'''|\leq |T''|$ and $|T^{'''}| = |T^{'''}\bigcap V(G)|+2$. It is straightforward to see that $T^{'''}$ is a $k$-strong shortest path union cover of  $G'$ and $T^{'''} \bigcap V(G)$ is a vertex cover of $G$.\\
\begin{figure}[h!] 
	\centering
	\includegraphics[scale=0.32]{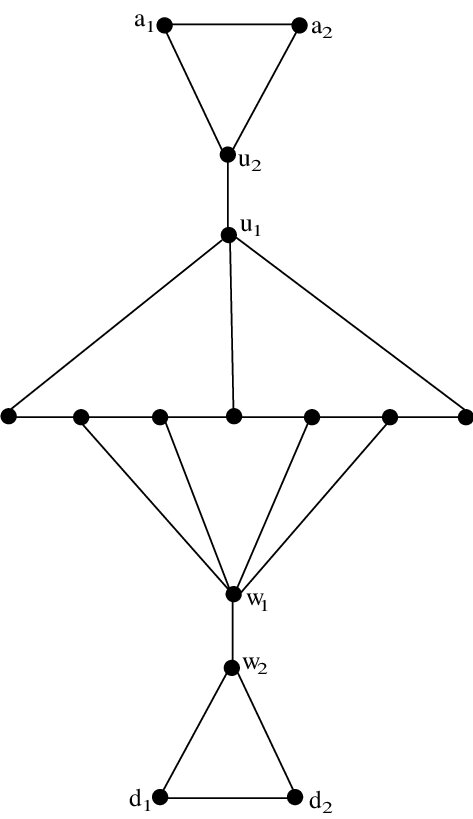}
	\caption {NP-complete illustration for a general graph}
	\label{fig:knpc3}
\end{figure}
		Case (ii): When $k=2$\\
		For $k=2$, Given a graph $G$ = $(V, E)$.  We construct a graph $G' = (V', E')$  from $G$ by subdividing every edge in $E(G)$ by a path of length $3$.  All the new vertices that are created in this step are called   path vertices. Now create  two new $K_{3}$  and introduce two new vertices $u_{1}$ and $w_{1}$ such that all the vertices of  $V(G)$ in $G^{'}$ will be joined with $u_{1}$ and all the path vertices in $G^{'}$ will be joined with vertex $w_{1}$. Join any one vertex from each $K_{3}$ to $u_{1}$  and $w_{1}$  as shown in the Figure 	\ref{fig:knpc3}. Thus the  graph $G'$ is obtained.\\
		\hspace*{1.0cm}  Let $S$ be the vertex cover of $G$. Then $S\cup \{a_{1},u_{1},d_{1}, w_{1}\}$ is a $2$-strong shortest path union cover.
		Conversely let $T$ be the $2$-strong shortest path union  cover. Let 
		$T^{'} = T\setminus \{a_{1}, a_{2}, u_{1}, u_{2},w_{1}, w_{2}, d_{1}, d_{2}\}  \cup \{a_{1}, u_{1}, d_{1}, w_{1}\}$. Clearly $T^{'}$ is a $2$-strong shortest path union cover and $|T'|\leq |T|$.\\ 
			If  $T'$ contains  a path vertex  $x'$, then replace $x'$ with 
	adjacent vertex $\textit{x} \in V$. Repeat this replacement procedure  until there are no more path vertices and the resulting set  be  $T^{''}$. Also $|T''|\leq |T'|$ and $|T^{''}| = |T^{''}\bigcap V(G)|+4$. It is  straightforward to see that 
	$T^{''}$ is a $2$-strong shortest path union cover of  $G'$ and $|T^{''}\bigcap V(G)|$ is a vertex cover of $G$.
		\end{proof}		
\end{thm}
\begin{rem}
The construction for $k=3,k=4$ and $k=5$ are illustrated in Figure \ref{fig:knpc2}.
	\end{rem}

\section{$2$-Strong shortest path union cover for certain graphs }
The following observations are easily verifiable.
\begin{obsr}
	For any path  $P_{m}$, $SSPC_{2U} (P_{m})$ = $\lceil \dfrac{ m}{5} \rceil$.
	\begin{proof}
		Let $P_m$  be the Path with order $m$. Because the maximam  degree of every vertex on $P_m$ is  $2$,  then degree $2$ vertex can cover at most $4$ edges at  distance $2$. Therefore the minimam number of vertices that can be coverd by $m$ vertices are $\lceil \dfrac{ m}{5} \rceil$. Thus,  $SSPC_{2U} (P_{m})$ $\geq$ $\lceil \dfrac{ m}{5} \rceil$. Also it can be easily verified that the set $S= \{v_{3},v_{7},\dots, v_{m-1}\}$ form the strong shortest path union cover set where $|S| = \lceil \dfrac{ m}{5} \rceil$. Therefore $SSPC_{2U} (P_{m})$ = $\lceil \dfrac{ m}{5} \rceil$.	
	\end{proof}		
\end{obsr}
\begin{obsr}
	For any cycle $C_{n}$, $SSPC_{2U} (C_{n})$ = $\lceil \dfrac{ n}{5} \rceil$, $n\geq5$.
	\begin{proof}
		Let $C_n$  be the Cycle with order $n$. Cycle is a  $2$ regular  graph. 
		Because the maximam  degree of every vertex on $C_n$ is  $2$,  then the  $2$ degree vertex can cover at most $4$ edges at  distance $2$. Therefore the minimam number of vertices that can be coverd by $n$ vertices are $\lceil \dfrac{ n}{5} \rceil$. Thus,  $SSPC_{2U} (C_{n})$ $\geq$ $\lceil \dfrac{ n}{5} \rceil$. Also it can be easily verified that the set $S= \{v_{1},v_{5},\dots, v_{n-3}\}$ form the strong shortest path union cover set where $|S| = \lceil \dfrac{ n}{5} \rceil$. Therefore $SSPC_{2U} (C_{n})$ = $\lceil \dfrac{ n}{5} \rceil$.	
	\end{proof}		
\end{obsr}
\begin{obsr}
	Let $A(m,n)$  be the Actinia graph for $m$ $\geq$ $2$ and $n$ $\geq$ $1$, then the $SSPC_{2U}(A(m,n))$ = $\lceil \dfrac{ n}{5} \rceil$, when $m$ is odd or even.
\end{obsr}
\begin{lem}
	For any Complete  bipartite graph  $K_{m, n}$,   $\ 2  \leq m \leq n$ , $SSPC_{2U} (K_{m, n})$ =  $m$.
	\begin{proof}
 Let $G$= $(M,N)$ be the complete bipartite graph $K_{m,n}$. 
 Let $2\leq m\leq n$. Let $U$ = $\{u_{1},u_{2},\dots,u_{m}\}$ and $W$ = $\{w_{1},w_{2},\dots,w_{n}\}$ be a bipartition of the bipartite graph  $k_{m,n}$. $\lvert U \lvert$ = $m$, $\lvert W \lvert$ = $n$.  Assume that  $T$ $\subseteq$$V(G)$ is the  $2$-strong shortest path union cover of $G$ and  $\lvert T\lvert$ = $m-1$.\\
case(i):  $\lvert T \lvert$ $\subseteq$ $U$.\\ 
Without loss of generality, $u_{m}$$\notin$ $T$. Let $T$ = $\{u_{1},u_{2},\dots,u_{m-1}\}$. Clearly the path between the vertices of $T$ and $u_{m}$ can cover atmost $(m-1)$ edges at distance $2$ adjacent with $u_{m}$. These paths will not cover $n-(m-1)$ edges incident with $u_{m}$. Also it is straightforward to see that these $n-(m-1)$ edges  incident with $u_{m}$ are not covered at distance $2$ by any other shortest path between the vertices of  $T$ and $N$.
Therefore $T$ is not the strong shortest path union cover of $K_{m,n}$.\\
case(ii) : $ T$ $\subseteq$ $W$. \\
Case(iii): $T$ = $T_{1}$ $\cup$ $T_{2}$.\\
The proof of  Case (ii) and (iii) are similar to case (i).\\ Again $T$ is not the $2$-strong shortest path union cover of $G$ by a similar argument.
Clearly the vertices of $M$  form the strong shortest path union cover. Hence $SSPC_{U}(K_{m,n})$ = $|M|$.				
	\end{proof}	
\end{lem}

\begin{figure}[h!]
	\centering 	
	\includegraphics[scale=0.5]{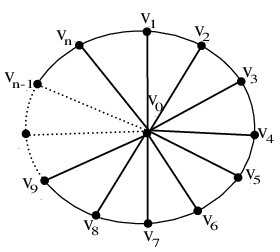}
	\caption{Wheel Graph}
	\label{fig:Wg}
\end{figure}

\begin{prop}\label{wh}
	For any Wheel graph  $W_{n}$, $n\geq5$, $SPC_{2U} (W_{n})$ = $SSPC_{2U} (W_{n})$ = $\lceil \dfrac{ n}{5} \rceil$.  
	\begin{proof}
		Let $v_{0}$ be the core vertex of  $W_{n}$. Since $v_{0}$ is adjacent to all the vertices and the diameter of the wheel graph is $2$ for $n \geq 5$.
	For $n$$\geq$ $2$ the graph $K_{1}$+$C_{n}$ is called wheel graph of order $n+1$ with $2n$ edges. The vertex on $K_{1}$ is called the core vertex denoted by $v_{0}$. The vertices that lie on the cycle $C_{n}$ is denoted by $v_{1},v_{2}, \dots, v_{n}$ in clockwise. Label the interier edges in anticlockwise $v_{0}v_{1}$, $v_{0}v_{n}$, $v_{0}v_{n-1}$, $\dots$, $v_{0}v_{2}$ by $1$, $3$, $5$, $\dots$, $2n$-$1$. Label the cycle edges in clockwise  $v_{1}v_{2}$, $v_{2}v_{3}$, $v_{3}v_{4}$, $\dots$, $v_{n}v_{1}$ by $1$, $2$, $3$, $\dots$, $n$. The $2$-strong shortest path union cover problem is to cover every edge of a graph at distance $2$ by the  unique shortest path from a  subset of vertices in the graph. Let S$\subseteq$ $V(W_{n})$ be the $2$-strong shortest path union cover of $W_{n}$. $S$ = \{$v_{1}$, $v_{5}$, $v_{9}$, $v_{13}$, $\dots$, $v_{n-3}$\}. Consider a set $T$ = $\{v_{i}, v_{i+1}, v_{i+2}, v_{i+3}, v_{i+4}\}$$\subseteq$ $V(G)$, $1\leq i\leq n$, then the vertex $v_{i+2}$ $\subseteq$ $S$ should be in $S$ such that the edges connected  $v_{i+2}$ to the remaining vertices in $T$/ \{$v_{i+2}$\} are covered in a unique path at distance  $2$. Continuing this for every distinct set of $5$ vertices in $C_{n}$ which is $V(W_{n}) /\{v_{0}\}$, such that one vertex ought to be chosen to cover all the edges at distance $2$. The  $SSPC_{2U} (W_{n})$  set of graph $W_{n}$, $n\geq5$ is a set with vertices $\lceil \dfrac{ n}{5} \rceil$ from the set $V(W_{n}) /\{v_{0}\}$.  Hence the  $2$-strong shortest path union cover for Wheel graph $W_{1,n}$ is $\lceil$$\dfrac{ n}{5} \rceil$.		
	 \end{proof}	
 \end{prop}

\begin{figure}[H] 
	\centering
	\includegraphics[scale=0.4]{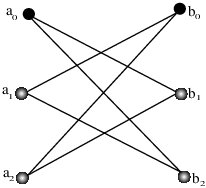}
	\caption {The vertices $\{a_0,b_0\}$ form  $SSPC_{2U} (H_{3,3})$}
	\label{fig:CG}
\end{figure}

\begin{prop}
For any Double Wheel graph  $DW_{n}$, $n\geq5$, $SPC_{2U} (DW_{n})$ = $SSPC_{2U} (DW_{n})$ = $2 \lceil \dfrac  { n}{5} \rceil$.  
\begin{proof}
	Follows from proposition \ref{wh}
 \end{proof}	
 \end{prop}
\begin{prop}
	For any Crown graph  $H_{n, n}$, $n \geq3$, $SPC_{2U} (H_{n, n}) = 2$.  
	\begin{proof}
		Let $H_{n, n}$, $n \geq3$ be a crown graph with vertex set $\{(a_i,b_j);0 \leq i,j\leq n-1, i \neq j\}$ of degree $n-1$. To cover all the edges at distance $2$ in shortest path require $2$ vertices. Figure \ref{fig:CG}, shows that the vertices $\{a_0,b_0\}$ form the $2$-strong shortest path union cover for $H_{n, n}$. If the $2$-strong shortest path union cover set  contains only any one vertex, there exists $n$ edges that are left uncovered. Thus $SPC_{2U} (H_{n, n}) = 2$.
	 \end{proof}	
\end{prop}
\begin{lem}
	For any Crown graph  $H_{n, n}$, $n \geq3$, $SSPC_{2U} (H_{n, n}) = n-1$  
	\begin{proof}
		Let $G$ be a crown graph $H_{n,n}$, where $n$ $\geq$ $3$. Let $U$ = $\{u_{1},u_{2},\dots,u_{n}\}$ and $V$ = $\{v_{1},v_{2},\dots,v_{n}\}$ be the vertex sets of a crown graph  $H_{n,n}$.
		Assume that  $T \subseteq V(G)$ is the  $2$-strong shortest path union cover of $G$ and  $\lvert T\lvert$ = $n -2$.\\
		Case(a):  $\lvert T \lvert$ $\subseteq$ $U$.\\
		 Without loss of generality, $u_{n}$$\notin$ $T$. Let $T$ = $\{u_{1},u_{2},\dots, u_{n-2}\}$. Clearly the path at distance $2$ between the vertices of $T$ and $u_{n}$ cover atmost  ${n-2}$ edges adjacent with $u_{n}$. The path between  $T$ and $u_{n}$ will not cover $2$ edges incident with  $u_{n}$. Also it is straightforward to see that these $2$ edges incident with $u_{n}$ are not covered at distance $2$ by any other strong shortest path between $T$ and $V$.
		Therefore $T$ is not a strong shortest path union cover of $H_{n,n}$.\\
		Case(b) : $ T$ $\subseteq$ $V$. \\
		Case(c): $T$ = $T_{1}$ $\cup$ $T_{2}$.\\ The proof of  Case (b) and (c) are similar to case (a).\\
		Clearly  $(n-1)$ vertices form  the  $2$-strong shortest path union cover  for  $H_{n,n}$.  Hence the proof.
	 \end{proof}	
\end{lem}

\begin{lem}
Let $G(n,k)$ be a  Petersen graph, then $SSPC_{2U} (G_{n, k}) = 3$.
	\begin{proof}
Let $G(n,k)$ be a Petersen graph. Let $a$, $b$, $c$ $d$ and $e$  be the outer vertices and $f$, $g$,  $h$, $i$ and $j$  be the corresponding inner vertices. Refer Figure \ref{fig:PGR}.\\
Assume that $T\subseteq V(G)$ is the $2$-strong shortest path union cover of $G$, and $|T| = 2$.\\
Case (i): $T$ contains any one inner and outer vertex\\
When $a \in T$, the edge $cd$ is left uncovered. To cover the edge $cd$  either $h$ or $i \in T$. Without loss of generality, $h \in T$, then the edge $cd$ is uncovered by $a$. If $h \in T$, then $h$ will not cover the edge $gi$. This implies that $T$ is not the $2$-strong shortest path union cover for $G(n,k)$.\\
\begin{figure}[H] 
	\centering
	\includegraphics[scale=0.5]{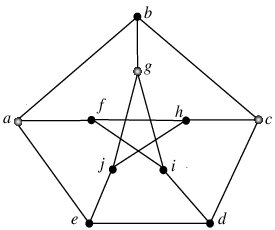}
	\caption {$S$ forms  $2$-strong shortest path union cover for Petersen graph}
	\label{fig:PGR}
\end{figure}
Case (ii): $T$ contains any two inner vertices.\\
Let $T = \{g,f \}$ where $g$ and $f$ be any two vertices from inner cycle.
There exists at least one edge $cd$ from outer cycle is left uncovered by $T$.
 This implies that $T$ is not the $2$-strong shortest path union cover for $G(n,k)$.\\
Case (iii): $T$ contains any two outer vertices\\
The proof is similar to the case (ii).\\
In these three cases, $T$ is not the $2$-strong shortest path union cover for $G(n,k)$.
Clearly $S= \{a,c,g \}$ forms the $2$-strong shortest path union cover for $G(n,k)$. And moreover  by Remark \ref{re1},  \ref{re3}, $SSPC_{U}G(n,k)$ = $3$. Hence it implies $SSPC_{kU}G(n,k)$ = $3$. 
\end{proof}	
\end{lem}

\begin{lem}
	For any Fan graph  $F_{1, n}$, $SPC_{2U} (F_{1, n})$ = $SSPC_{2U} (F_{1, n})$ = $\lceil \dfrac{ n}{5} \rceil$, $n\geq5$.
	\begin{proof}
	A fan graph $F_{1,n}$ is defined as the graph $K_{1}$+$P_{n}$, where $K_{1}$ is with single vertex and $P_{n}$ is a path on n vertices. Let $v_{o}$ be the central vertex  and $v_{1}$, $v_{2}$, $\dots$, $v_{n}$ be the vertices in  $P_{n}$.
	Let $S$ $\subseteq$ $V(G)$ be the $2$-strong shortest path union cover of $G$. Then $S$ contains any one vertex out of every $5$ vertices from $P_{n}$ which form the $2$-strong shortest path union cover of $G$. It is straightforward to see any $\lceil \dfrac{ n}{5} \rceil-1$ vertices will not cover $E(F_{1,n})$.
	Hence $SSPC_{2U} (F_{ 1,n})$ = $\lceil \dfrac{ n}{5} \rceil$. 	
	\end{proof}	
\end{lem}
\begin{lem}
	For any Double Fan graph  $DF_{ n}$,  $SPC_{2U} (DF_{ n}) =  SSPC_{2U} (DF_{ n}) =1+\lceil \dfrac{ n}{5} \rceil$  
	\begin{proof}
		The double fan graph  $F_{2,n}$ consists of $n+2$ vertices and $3n-1$ edges. The path vertices 
		 $v_{1}$, $v_{2}$, $\dots$, $v_{n}$ are adjacent to $v_{o}$ and $v_{o}'$. Let S$\subseteq$ $V(G)$ be the $2$-strong shortest path union cover of $G$. Then $S$ contains $v_{o}$ or $v_{o}'$  and any one vertex out of every $5$ vertices from $P_{n}$ which form a $2$-strong shortest path union cover of $G$.
		 Thus  $SPC_{2U} (DF_{ n}) =  SSPC_{2U} (DF_{ n}) =1+\lceil \dfrac{ n}{5} \rceil$.
	 \end{proof}	
\end{lem}

\begin{thm}
	For any Generalised  Friendship  graph  $F_{k,n}$, $k,n\geq4$,  $SSPC_{2U} (F_{k,n}) \leq 1+n \lceil\dfrac{k}{5}\rceil$
	\begin{proof}
	Let $v_{0}$ be the central vertex and  $v_{1},v_{2},\dots, v_{n}$ be the outer vertices of the petals. The central or core vertex $v_{0}$ covers the edges at distance $2$ in each petal. Still there exists edges within the petals left uncovered. By choosing a vertex out of every five vertices in each petal would cover the edges at distance $2$.  Hence there exists $1+n (\lceil\dfrac{k}{5}\rceil$) vertices for $2$-strong shortest path union cover in $F_{k,n}$. 		
	\end{proof}	
\end{thm}
\begin{rem}
	For any Generalised  Friendship  graph  $F_{k,n}$, $k,n\geq4$,  $SPC_{2U} (F_{k,n}) \leq 1+n \lceil\dfrac{k}{5}\rceil$.
\end{rem}

\section{$2$-Strong shortest path union cover for Networks}
\par \hspace*{0.5cm} Interconnection networks play a key role in the design and implementation of communication networks and the recent advent of optic technology add more design problems \cite{XJ01}. In general, an interconnection network may be modeled by a simple graph whose nodes represent components of the network and whose links represent physical communication links. The vertex which connects the two degree vertices in $(BF(2))$ known as a binding vertex.
In this section we determine the   $2$-shortest path and strong  shortest path union cover for various networks consist of  Butterfly network, Augmented Butterfly network, Enhanced butterfly network, Benes network, Silicate network $SL(n)$,  Hypercube $Q_{n}$, Sierpi\'{n}ski graphs and Sierpi\'{n}ski gasket graphs.

\begin{thm}
	Let $G$ be an $r$-dimensional butterfly network. Then 
	\begin{center}
		$SPC_{2U}(BF(r))$ $\leq$ 
		$\begin{cases} 2^{r-1}, when\ r\leq4 \\ 
		2^{r}, when\ r>4\\
		\end{cases}$
	\end{center}
\begin{figure}[H] 
	\centering
	\includegraphics[scale=0.15]{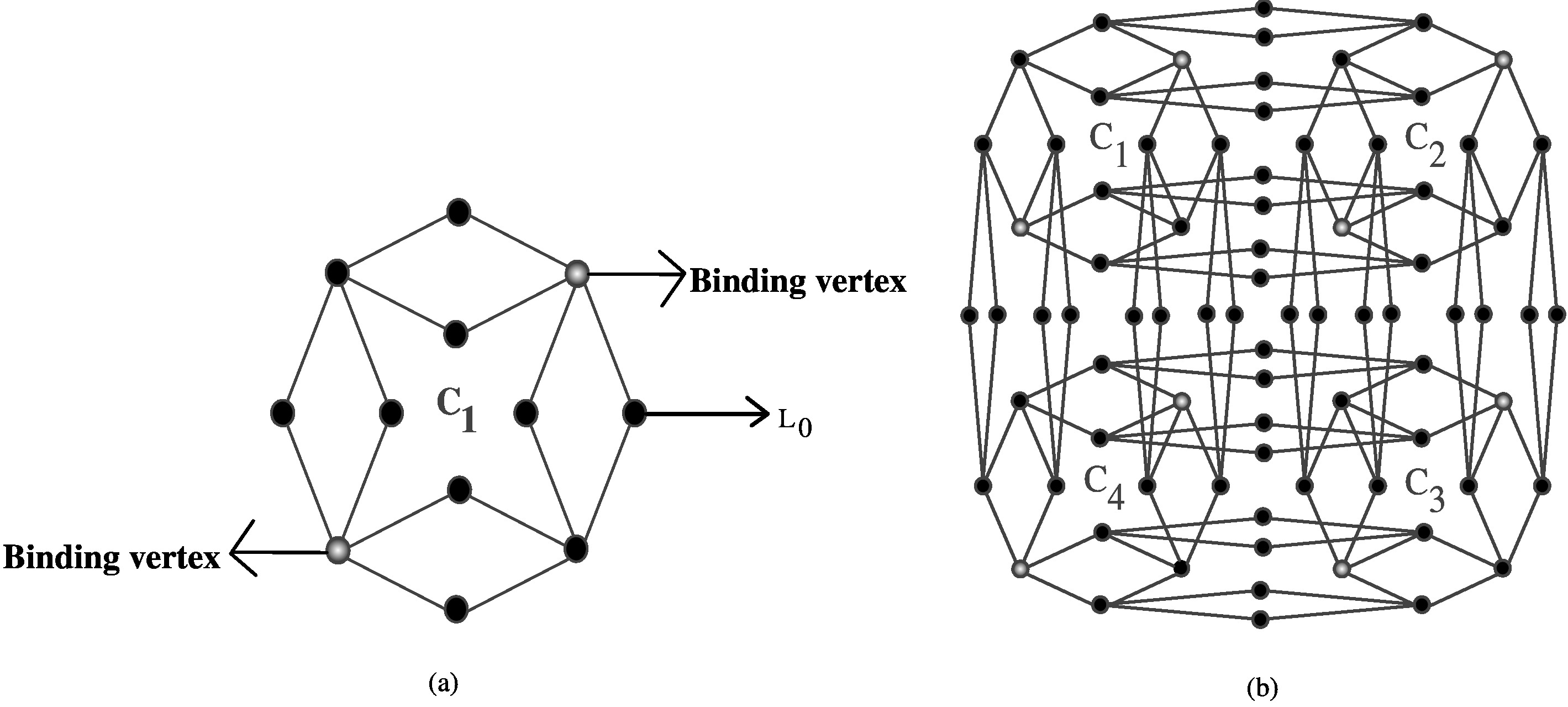}
	\caption {(a) Pairs of binding vertices form $2$-shortest path union cover for  $BF(2)$  (b) Pairs of binding vertices in each $BF(2)$ form $2$-shortest path union cover for $BF(4)$}
	\label{fig:BF}
\end{figure}
	\begin{proof}
		case(i) :  when $ r\leq4$\\
Consider butterfly network $BF(2)$. $BF(2)$ is a bipartite graph with $4$ vertices in each part, which are independent sets called as levels $L_0$ and $L_1$. The removal of level $0$ vertices of  $BF(2)$ leaves $2$ disjoint copies of  $BF(1)$. The binding vertex in each copy of  $BF(1)$ forms the $2$-shortest path union cover for $BF(2)$. Similarly for $BF(3)$, we have $2$ copies of $BF(2)$ where the binding vertex in each copy of $BF(2)$ forms the $2$-shortest path union cover for $BF(3)$. As the butterfly network has dual symmetry, in $BF(4)$ there are $2$ copies of  $BF(3)$ and $4$ copies of  $BF(2)$. The pairs of  binding vertices in each  $BF(2)$ form the $2$-shortest path union cover for $BF(4)$. Hence $SPC_{2U}(BF(4)) = 8$. Therefore $SPC_{2U}(BF(r)) \leq 2^{r-1}$.\\	case(ii) : When $r>4$\\
		Consider the butterfly graph $BF(5)$. By the recursive construction  $BF(5)$ has $2$ copies of $BF(4)$ and a new level $L_4$ with $2^{5}$ vertices. Let $S_1$ and $S_2$  denote $2$-shortest path union cover sets in the left and the right copy of  $BF(4)$ respectively. Since $SPC_{2U}(BF(4)) = 8$, select all $16$ binding vertices from left copy of $BF(5)$ into $S_1$ and $16$ binding vertices from right copy of $BF(5)$ into $S_2$. Let $S_1$= $\{(1,i)(2,i)/ $i=0 to 15$\}$. Then all the edges in left copy of $L_1$  are covered at distance $2$ by the binding vertices in $S_1$. Using the symmetry of butterfly graphs,  mirror image compliment set of $S_1$ in the right copy of  $BF(5)$ will cover all the edges at distance $2$ in that copy and it is given by $S_2$= $\{(1,i)(2,i)/ $i=16 to 31$\}$. Hence $S$= $S_1 \cup S_2$ form the $2$-shortest path union cover for $BF(5)$ with $\rvert S\rvert$ = $2.2^{4}$ = $2^{5}$. By symmetry, the result holds true for $n > 5$. Therefore $SPC_{2U}(BF(r)) \leq 2 ^{r}$, when $\ r \geq 5$.
	\end{proof}	
\end{thm}
\begin{figure}[h] 
	\centering
	\includegraphics[scale=0.5]{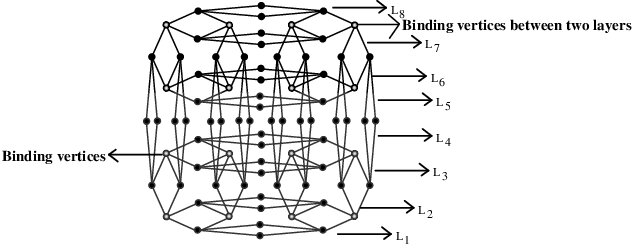}
	\caption {Binding vertices in each $BF(2)$ form $2$- strong shortest path union cover for $BF(4)$}
	\label{fig:BF1}
\end{figure}
\begin{thm}
	Let $G$ be an $r$-dimensional butterfly network. Then  $SSPC_{2U}(BF(r)) \leq \lceil \dfrac{r}{2}\rceil 2^{r-1}$, $n\geq3$.
	\begin{proof}
		Let $G$ be an $r$-dimensional butterfly network. By recursive construction, $BF(r)$ has $2$ copies of  $BF(r-1)$. Let $S_{1}$ and $S_{2}$ denote the $2$-strong shortest path union cover sets in the bottem and top copy of $B(r-1)$ respectively. In $BF(3)$ there exists $2$ layers of bottem and top layers. The binding vertices between two consecutive layers in bottom and top  cover all the edges at distance $2$ and form the $2$-strong shortest path union cover for $BF(3)$. Similarly in $BF(4)$, there exists four  bottom and top layers. The binding vertices between two consecutive layers both in bottom and top cover all the edges at distance two and form the $2$-strong shortest path union cover for $BF(4)$. Therefore $SSPC_{2U}(BF(4)) = 2.2^{3} = \lceil \dfrac{r}{2} \rceil 2^{r-1}$.\\
		Proceeding like this, for $BF(r)$, there exists bottem  and top copy of $BF(r-1)$. The binding vertices  between two consecutive layers cover all the edges at distance $2$ in bottem and top copies of  $BF(r-1)$. Hence $S=S_{1}\cup S_{2}$ becomes the $2$-strong shortest path union cover of  $BF(r)$ with $\rvert S \rvert \leq \lceil \dfrac{r}{2} \rceil 2^{r-1}$. Therefore $SSPC_{2U}(BF(r)) \leq  \lceil \dfrac{r}{2} \rceil 2^{r-1}$.		
	\end{proof}	
\end{thm}

\begin{thm}\label{ag}
	Let $G$ be an $3$-dimensional Augmented butterfly network. Then  $SSPC_{2U}(ABF(3)) =  12$.

	\begin{figure}[h] 
		\centering
		\includegraphics[scale=0.5]{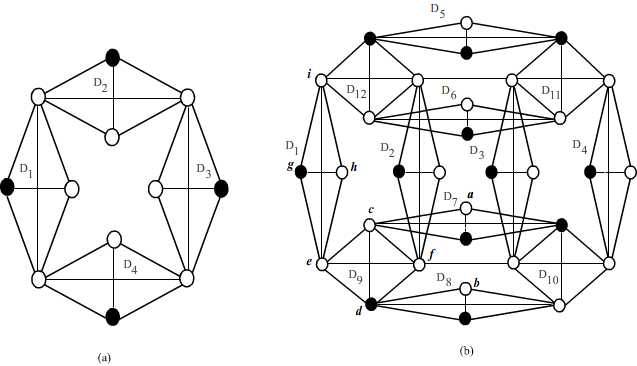}
		\caption {(a) and (b) One vertex from each diamond  form  $2$-strong shortest path union cover for $ABF(2)$ and $ABF(3)$}. 
		\label{fig:ABN}
	\end{figure}
	\begin{proof}	
		Consider the $r$-dimensional butterfly network $BF(r)$. Place an new edge on the antipodal vertices in a cycle. The resulting graph is known as Augmented butterfly network.\\ Let $G$ be	an $3$-dimensional Augmented butterfly network consists of diamond representation as shown in Figure \ref{fig:ABN}. Let $S$ be the $2$-strong shortest path union cover for  $ABF(3)$. Let $a,b\in S$. The vertices $a,b$ can cover all the edges of the diamond $cedf$, except  the edge  $ef$. Therefore to cover the edge $ef$  we require  at least one vertex from the diamond $cedf$. Thus for each diamond of $ABF(3)$ we have to take at least one vertex in the $2$-strong shortest path union cover set. 	Hence $SSPC_{2U}(ABF(3)) = 12$.
	\end{proof}	
\end{thm}
	\begin{cor}
		Let $G$ be an $r$-dimensional Augmented butterfly network. Then the $SSPC_{2U}(ABF(r)) \leq r (2^{r-1})$.
		\begin{proof}	
		By Theorem \ref{ag}, for $SSPC_{2U}(ABF(3))  = 12$. Similarly for $SSPC_{2U}(ABF(4))  = 28$.  Proceeding like this for $r$-dimension of augmented butterfly network,  there exists $ r (2^{r-1})$ diamond in $(ABF(r))$. Also it is straightforward to note that a set $S\subset V(G)$ that contains one vertex from each diamond of $ABF(r)$ forms the $2$-strong shortest path union cover. 
		Hence $SSPC_{2U}(ABF(r)) \leq r (2^{r-1})$.
	\end{proof}	
\end{cor}
\begin{rem}
Let $G$ be an $r$-dimensional augmented butterfly network. Then  $SPC_{2U}(ABF(r)) \leq r (2^{r-1})$.
\end{rem}
\begin{thm}\label{ebg}
	Let $G$ be an $3$-dimensional Enhanced butterfly network. Then  $SSPC_{2U}(EBF(3)) = 12$.
	\begin{figure}[h] 
		\centering
		\includegraphics[scale=0.5]{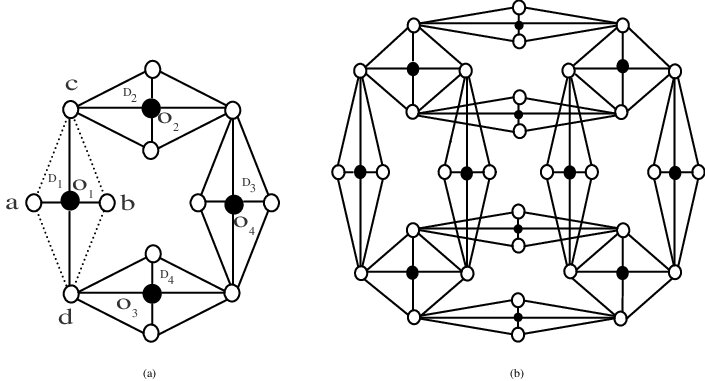}
		\caption {(a) and (b) The core vertices form $2$-strong shortest path union cover for $EBF(2)$ and $EBF(3)$}
		\label{fig:EBN}
	\end{figure}
	\begin{proof}
		Consider the $r$-dimensional butterfly network $BF(r)$. Place a new vertex in each $4$-cycle of $BF(r)$ and join this vertex to the four vertices of the $4$-cycle. The resulting graph is called an enhanced butterfly network $EBF(r)$.\\
		Let $G$ be an $3$-dimensional enhanced butterfly network consists of diamond representation as shown in Figure  \ref{fig:EBN}. The center vertex in each diamond is known as core vertex $o$. The centre vertex with four degree in diamond structure $D_{1}: cbda$ known as core vertex $o_{1}$ as shown in Figure  \ref{fig:EBN}(a). Let $S$ be the $2$-strong shortest path union cover for  $EBF(2)$. Let $o_{2}, o_{3} \in S$. The core  vertices $o_{2}, o_{3}$ can cover all the edges of the diamond $cbda$, except  the edges $(ao_{1})$ and $(o_{1}b)$ in  $D_{1}$. These edges is not covered by any other vertices in the diamond $cbda$, except the core vertex $o_{1}$. Hence in each diamond structure of  $EBF(2)$,  core vertex should be there in $2$-strong shortest path union cover set. Therefore for $SSPC_{2U}(EBF(2)) = 4$. Similarly for $3$ dimensional enhanced butterfly network, we require a core vertex in each diamond. Hence $SSPC_{2U}(EBF(3)) = 12$. 
		\end{proof}	
\end{thm}	

	\begin{thm}
		Let $G$ be an $r$-dimensional Enhanced butterfly network. Then  $SSPC_{2U}(EBF(r)) \leq r (2^{r-1})$.
		\begin{proof}
		By Theorem \ref{ebg}, $SSPC_{2U}(EBF(3)) = 12$. Similarly for $SSPC_{2U}(EBF(4)) = 28$.
		Proceeding like this for $r$-dimension of enhanced butterfly network,  there exists $ r (2^{r-1})$ diamond in $EBF(r)$. Also it is straightforward to note that a set $S\subset V(G)$ that contains core vertex from each diamond of $EBF(r)$ forms the $2$-strong shortest path union cover. 
		Hence $SSPC_{2U}(EBF(r)) \leq r (2^{r-1})$.

	\end{proof}	
\end{thm}
\begin{rem}
	Let $G$ be an $r$-dimensional enhanced butterfly network. Then  $SPC_{2U}(EBF(r)) \leq r (2^{r-1})$.
\end{rem}	
\begin{cor}
	Let $G$ be an $r$-dimensional benes network. Then  $SPC_{2U}(B(r)) \leq 2^{r}$ for $r \geq 3$.
	\begin{proof}
		Let $G$ be an $r$-dimensional benes network. The removal of level $0$ vertices $v_{1}, v_{2}, \dots, v_{n}^{r}$ of $B(r)$ gives two subgraphs $G_{1}$ and $G_{2}$	of  $B(r)$, each isomorphic to $B(r-1)$. The pairs of binding vertices  in each $B(2)$ covers all the edges in $B(r)$. No other vertices can cover the edges in the layers except the binding vertices. Hence the pairs of binding vertices beween layers form $2$-shortest parh union cover for $B(r)$. 
		Hence $SPC_{2U}(B(r)) \leq 2^{r}$ for $r \geq 3$.
	\end{proof}	
\end{cor}
\begin{figure}[h] 
	\centering
	\includegraphics[scale=0.3]{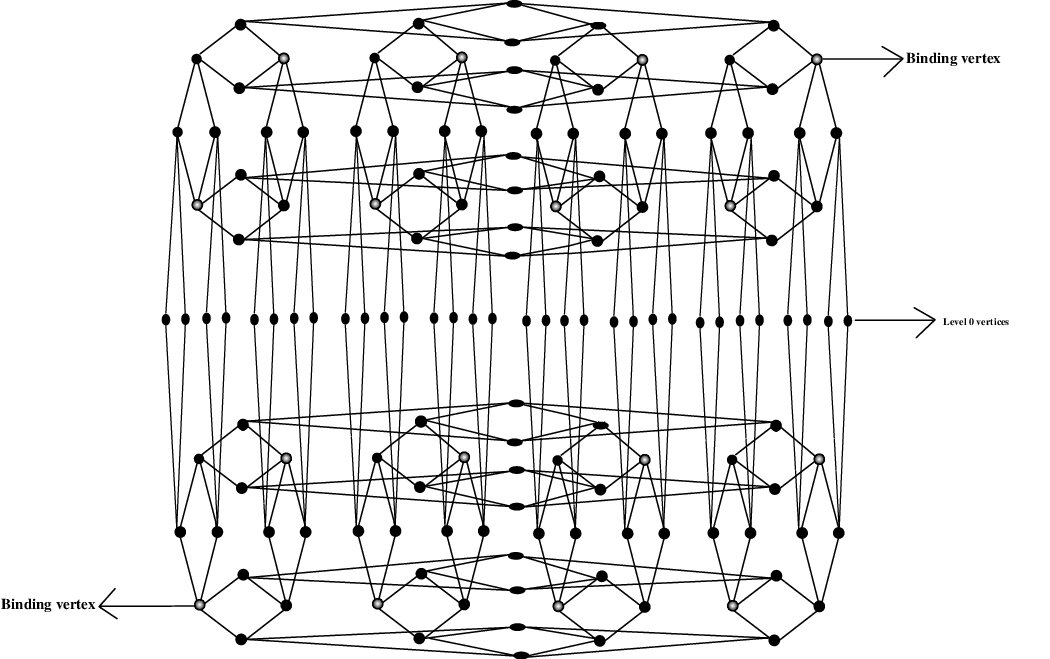}
	\caption {Binding vertices in each $B(2)$ form $2$- shortest path union cover for $B(4)$}
	\label{fig:BN}
\end{figure}

\begin{figure}[h] 
	\centering
	\includegraphics[scale=0.3]{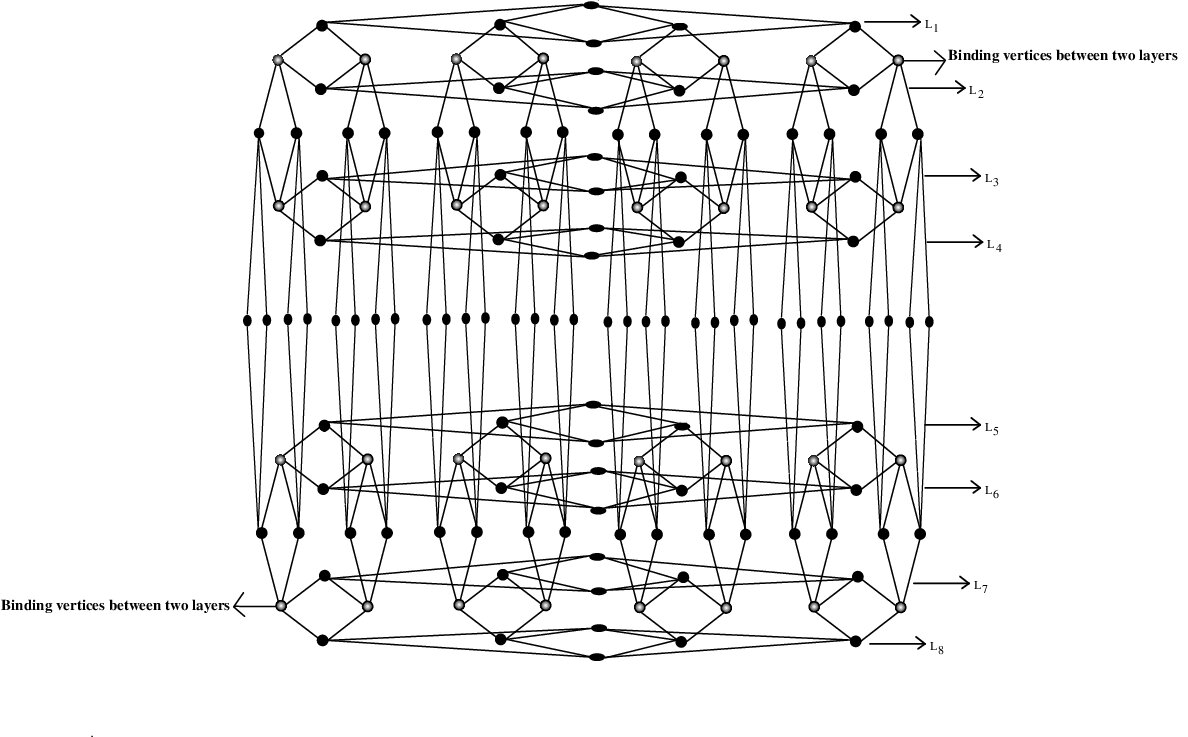}
	\caption {Binding vertices between two layers form $2$-strong shortest path union cover for  $B(4)$}
	\label{fig:BN3}
\end{figure}
\begin{rem}
	For an $2$-dimensional benes network, the $SPC_{2U}(B(r)) =2$.\\ 
	The two antipodal vertices  $u_{1}, u_{2}$ form the $2$-shortest path union cover for $B(2)$. 
\end{rem}

\begin{thm}
	Let $G$ be an $r$-dimensional benes network. Then  $SSPC_{2U}(B(r)) \leq  \lceil \dfrac{r}{2} \rceil 2^{r}$
	\begin{proof}
		The benes network consists of back to back butterflies. Let $G$ be an $r$-dimensional benes network. By recursive construction, $B(r)$ has $2$ copies of 	$B(r-1)$. Let $S_{1}$ and $S_{2}$ denote $2$-strong shortest path union cover sets in the bottem and top copy of $B(r-1)$ respectively. In $B(2)$ there exists $2$ layers of bottem and top layers. There are two pairs antipodal vertices in bottom and top layers which cover all the edges at distance $2$ and form the $2$-strong shortest path union cover for $B(2)$. Similarly in $B(3)$ and 
		 $B(4)$, the binding vertices  between two consecutive layers  in bottom and top copy of $B(3)$ cover all the edges at distance $2$, thus form the $2$-strong shortest path union cover. Therefore $SSPC_{2U}(B(4)) = 2.2^{4} = \lceil \dfrac{r}{2} \rceil 2^{r}$.\\
		Proceeding like this, for $B(r)$, there exists bottom and top copy of $B(r-1)$. The  binding vertices  between two consecutive layers cover all the  edges at distance $2$ in bottom and top copies of  $B(r-1)$. Hence $S=S_{1}\cup S_{2}$ becomes the $2$-strong shortest path union cover of  $B(r)$ with $\rvert S \rvert \leq \lceil \dfrac{r}{2} \rceil 2^{r}$. Therefore $SSPC_{2U}(B(r)) \leq  \lceil \dfrac{r}{2} \rceil 2^{r}$.	
	\end{proof}	
\end{thm}

\begin{figure}[h] 
	\centering
	\includegraphics[scale=0.5]{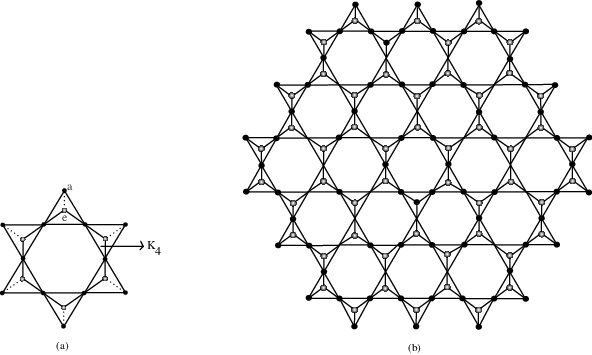}
	\caption {(a)  and (b) Three degree vertices in each $K_4$ form the $2$-shortest and strong shortest path union cover  for  $SL(1)$ and   $SL(3)$ }
	\label{fig:SN}
\end{figure}

\begin{thm}\label{5.3}
	Let $G$ be the silicate network $SL(n)$  of dimension $n$, then  $SSPC_{2U}(G) = 6n^{2}$
	\begin{proof}
	Let $SL(n)$  be the silicate network with $15n^{2}$+$3n$ vertices and $36n^{2}$ edges. Let $S\subseteq V(SL(n))$ be the $2$-strong shortest path union cover set for $SL(n)$. In each $K_{4}$ of $SL(1)$, there exists an edge $(ae)$ which is uncovered by any other vertices in $SL(1)$. Hence the $3$ degree simplicial vertices in each $K_{4}$ should be in $2$-strong shortest path union cover set. Therefore $S$ contains only $3$ degree simplicial vertices in each $K_{4}$ of $SL(n)$ as shown in Figure 	\ref{fig:SN}(a). 
	 Therefore $S$  contains the  $6n^{2}$ vertices  in $SL(n)$ and form the $2$-strong shortest path union cover set for $SL(n)$. Refer Figure \ref{fig:SN}(b). Hence  by Theorem \ref{cli}, $SSPC_{2U}(G) = 6n^{2}$.	
	\end{proof}	
\end{thm}
\begin{rem}
	Let $G$ be the silicate network $SL(n)$  of dimension $n$, then  $SPC_{2U}(G) = 6n^{2}$.
	\end{rem}
\begin{figure}[h] 
	\centering
	\includegraphics[scale=0.5]{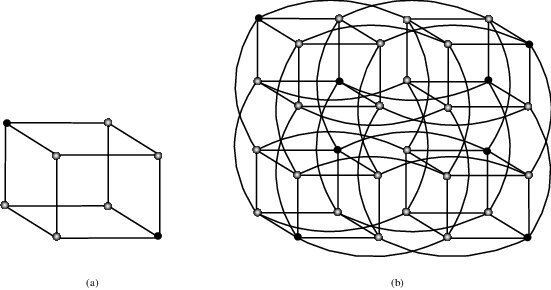}
	\caption {(a) and (b)The set of  vertices marked in dark form the $2$-strong shortest path union cover for $Q_{3}$ and  $Q_{5}$}
	\label{fig:HYN}
\end{figure}
\begin{thm}\label{5.4}
	For the $n$ cube $Q_{n}$,  $n\geq3$, then  $SSPC_{2U}(Q_{n}) \leq 2^{n-2}$.
	\begin{proof}
	Let $Q_{n}$ be the hypercube network with $2^{n}$ verices. For $n=3$, 	$Q_{3}$ has two copies of $Q_{2}$. It can be easily verified that in each copy of $Q_{2}$ atleast one vertex is sufficiently enough to cover all the edges at distance $2$ to form the $2$-strong shortest path union cover in $Q_{3}$. For $Q_{4}$, there exists four copies of $Q_{2}$, one vertex in each copy of $Q_{2}$ is more than enough to cover all the edges at distance $2$ in $Q_{4}$. For $Q_{5}$, there exists eight copies of $Q_{2}$, atleast one vertex in each copy of $Q_{2}$ is more than enough to cover all the edges at distance $2$ in $Q_{5}$ as shown in Figure. Similarily proceeding by induction method for $Q_{n}$, there exists $2^{n-2}$ copies of $Q_{2}$. Choosing atleast one vertex in each copy of $Q_{2}$ cover all the edges at distance $2$ and form the $2$-strong shortest path union cover for $Q_{n}$. Thus it can  easily be verified that $2^{n-2}$ vertices are sufficient to form the $2$-strong shortest path union cover for $Q_{n}$. Hence by Remark \ref{re1},  \ref{re3},  
	$SSPC_{2U}(Q_{n})  \leq 2^{n-2}$.	
	\end{proof}	
\end{thm}
\begin{rem}\label{Re}
	For the $n$ cube $Q_{n}$,  $n\geq3$, then  $SPC_{2U}(Q_{n}) \leq 2^{n-2}$.
\end{rem}
\begin{rem}
The result in Theorem \ref{5.4} and Remark \ref{Re} is sharp for $n=3$.
\end{rem}
\begin{thm}
	Let $G$ be the Sierpi\'{n}ski Graph $S(n,3)$, $n$ $\geq$  $2$, then  $SSPC_{2U}(G) \leq 3^{n-1}$.
	
	\begin{figure}[h] 
		\centering
		\includegraphics[scale=0.4]{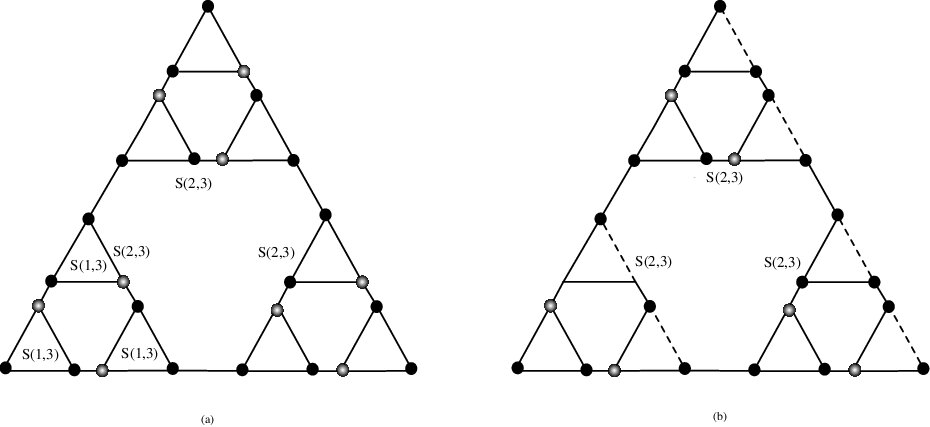}
		\caption {(a) $S$ forms $2$-shortest and strong shortest path union cover  for  $S(3,3)$ (b) The dotted lines represent the edges that are not covered  at distance $2$ by $T$ in  $S(3,3)$}.
		\label{fig:b}
	\end{figure}
	\begin{proof}
	Let $G(V,E)$ be the $n$-dimensional 	Sierpi\'{n}ski Graph $S(n,3)$, $n$ $\geq$ $2$. There exists  $3^{(n-1)}$ copies of $S(1,3)$, in $S(n,3)$. 
	In $S(2,3)$, the three copies of  $S(1,3)$	are connected to each other by an edge as depicted in Figure \ref{fig:b} (a). In $S(3,3)$, the three copies of $S(2,3)$ are connected to each other by an edge. Similarly in $S(n,3)$, the three copies of $S(n-1,3)$ are connected by an edge. In each  $S(2,3)$, choose the alternative   vertices of $C_{6}$ except the extreme  vertices in $S(2,3)$ which cover the edges at distance $2$ and  form a $2$-strong shortest path union cover as shown in \ref{fig:b}(a). The similar pattern of choosing the vertices in  each  $S(2,3)$ is followed in $S(n,3)$. 
	Proceeding like this, there exists $3^{(n-2)}$ copies of  $S(2,3)$ in $S(n,3)$. Let $S$  be the $2$-strong shortest path union cover of $G$ and $S$ contains $3^{(n-1)}$  vertices which cover the edges at distance $2$ in $S(n,3)$.\\
	Assume that $T$ $\subseteq$ $V(G)$
	such that $|T|$ $\leq$ $|S|$ forms the $2$-shortest and strong shortest path union cover.
	${T}$ contains two alternative vertices of   $C_{6}$  in $S(2,3)$.\\
	There exists $2$ edges left uncovered in each copy of $S(2,3)$ in $S(n,3)$ by the vertices in $T$ as shown in \ref{fig:b} (b). This implies that $T$ does not form $2$-strong shortest path union cover for $S(n,3)$. Hence $S$ forms a $2$-strong shortest path union cover and by Remark \ref{re1},  \ref{re3}, $ SSPC_{2U}(G) \leq 3^{n-1}$ 
\end{proof}	
\end{thm}
\begin{rem}
Let $G$ be the Sierpi\'{n}ski Graph $S(n,3)$, $n$ $\geq$  $2$, then  $SPC_{2U}(G) \leq 3^{n-1}$.
\end{rem}

\begin{figure}[h] 
	\centering
	\includegraphics[scale=0.4]{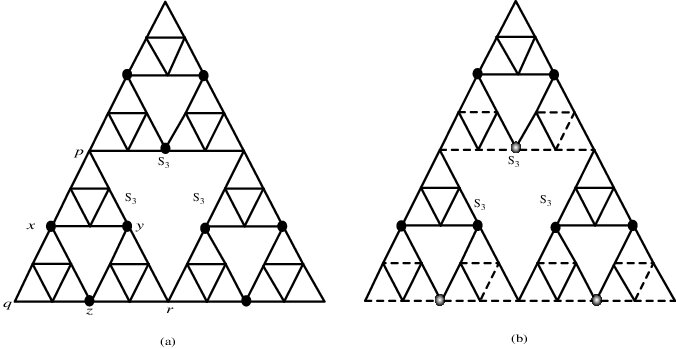}
	\caption {(a) $S$ forms the $2$-shortest path union cover  for  $S_{4}$ (b) The dotted lines represent the edges that are not covered  at distance $2$ by $T$ in  $S_{4}$}.
	\label{fig:d}
\end{figure}
\begin{thm}
	Let $G$ be the Sierpi\'{n}ski gasket Graph $S_{n}$, $n$ $\geq$  $3$, then  $SPC_{2U}(S_{n}) \leq 3^{n-2}$.
	\begin{proof}
		Let $S_{n}$ be the $n$- dimensional Sierpi\'{n}ski  gasket Graph $S_{n}$, $n$ $\geq$  $3$. 
		In $S_{3}$, the three copies of  $S_{2}$ are merged to each other. In $S_{4}$, the three copies of $S_{3}$, are merged to each other as depicted in Figure \ref{fig:d} (a). Similarly in $S_{n}$, the three copies of $S_{n-1}$, are merged to each other. In $S_{3}$,  choosing any three vertices except the merged vertices would leave at least one edge in each $K_{3}$ uncovered. Hence choose three merged vertices $\{x,y,z \}$ which cover the edges at distance $2$ and  form the $2$-shortest path union cover as shown in \ref{fig:d}(a). The similar pattern of choosing the vertices in  each  $S_{3}$  is followed in $S_{n}$. 
		Proceeding like this, there exists $3^{(n-3)}$ copies of $S_{3}$, in $S_{n}$. Let $S$  be the $2$-shortest path union cover for $G$ and $S$ contains three merged vertices in each copy of $S_{3}$ which cover the edges at distance $2$ in $S_{n}$. Hence  $3^{n-2}$ vertices in $S$ form the $2$-shortest path union cover for $S_{n}$.\\
		Assume that $T$ $\subseteq$ $V(G)$
		such that $|T|$ $\leq$ $|S|$ forms the $2$-shortest  path union cover.
		 ${T}$ contains any two merged vertices in  each $S_{3}$ of  $S_{n}$.\\
		There exists $7$ edges left uncovered in each  $S_{3}$  of $S_{n}$ by the vertices in $T$ as shown in \ref{fig:d}(b) . This implies that $T$ does not form the $2$-shortest path union cover for $S_{n}$. Hence $S$ forms the $2$-shortest path union cover for $S_{n}$  and by Remark \ref{re1},  \ref{re3}, $SPC_{2U}(S_{n}) \leq 3^{n-2}$.
			\end{proof}	
	\end{thm} 

\begin{figure}[h] 
	\centering
	\includegraphics[scale=0.4]{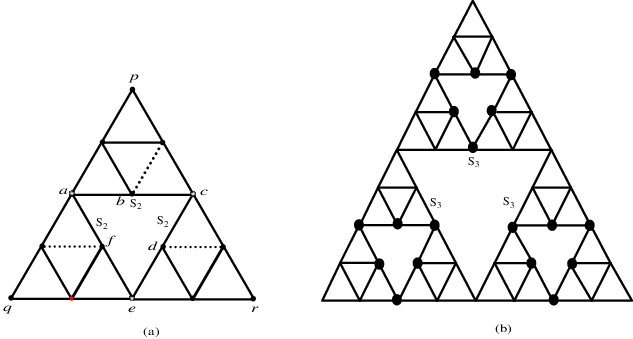}
	\caption {(a) The dotted lines represent the edges that are not covered  at distance $2$ by $T$ in $S_{3}$ (b)  $S$ forms $2$-strong shortest path union cover  for  $S_{4}$}.
	\label{fig:e}
\end{figure}
\begin{thm}
	Let $G$ be the Sierpi\'{n}ski gasket Graph $S_{n}$, $n$ $\geq$  $3$, then  $SSPC_{2U}(S_{n}) \leq 6(3^{n-3})$.
	\begin{proof}
		Let $S_{n}$ be an $n$- dimensional Sierpi\'{n}ski  gasket Graph $S_{n}$, $n$ $\geq$  $3$. 
		In $S_{3}$, the three copies of  $S_{2}$ are merged to each other and in  $S_{4}$, the three copies of $S_{3}$, are merged to each other. Similarly in $S_{n}$, the three copies of $S_{n-1}$, are merged to each other. In each $S_{3}$, choose three merged vertices and the vertices between the merged vertices  which cover the edges at distance $2$ and  form the $2$-strong shortest path union cover as shown in \ref{fig:e}(b). The similar pattern of choosing the vertices in  each copy of  $S_{3}$  is followed in $S_{n}$. 
		Proceeding like this, there exists $3^{(n-3)}$ copies of $S_{3}$, in $S_{n}$. Let $S$  be the $2$-strong shortest path union cover of $G$ and $S$ contains three merged vertices and the vertices between the merged vertices  in each copy of  $S_{3}$ which cover the edges at distance $2$ in $S_{n}$.\\
		Assume that $T$ $\subseteq$ $V(G)$
		such that $|T|$ $\leq$ $|S|$ forms the $2$ shortest and strong shortest path union cover.\\
		${T}$ contains only the merged vertices in each $S_{3}$ of  $S_{n}$.\\
		There exists at least one edge  left uncovered in each copy of  $S_{2}$ in $S_{n}$ by the vertices in $T$ as shown in \ref{fig:e}(a). This implies that $T$ does not form the $2$-strong shortest path union cover for $S_{n}$. Hence $S$ forms the $2$-strong shortest path union cover for $S_{n}$ and by Remark \ref{re1},  \ref{re3}, $SSPC_{2U}(S_{n}) \leq 6(3^{n-3})$.
	\end{proof}	
\end{thm}
\begin{rem}
	Let $G$ be the Sierpi\'{n}ski gasket Graph $S_{n}$, $n$ =  $2$, then $SPC_{2U}(S_{n}$)=$SSPC_{2U}(S_{n})=2$.
\end{rem}	
\section {Conclusion}
In this manuscript we have determined the complexity results for $k$-strong shortest path union cover, the $2$-strong shortest path union cover for general graphs, various networks,  Sierpi\'{n}ski graphs and Sierpi\'{n}ski gasket graphs.  Further $k$-strong shortest path union cover for other product graphs and networks are under consideration.

\end {document}